\documentclass[12pt] {article}
\usepackage{amsmath,amsthm, amscd}
\usepackage{amssymb,latexsym}

\title{Plurisubharmonic functions on hypercomplex manifolds and
HKT-geometry.}
\author{Semyon Alesker\footnote{Partially supported by ISF grant 1369/04.}, Misha Verbitsky\footnote{Misha Verbitsky is
an EPSRC advanced fellow supported by CRDF grant RM1-2354-MO02 and
EPSRC grant  GR/R77773/01} }

\date{}
\def\RR{\mathbb{R}}
\def\CC{\mathbb{C}}

\def\HH{\mathbb{H}}

\def\Ome{\Omega}
\def\ome{\omega}
\def\lam{\lambda}
\def\Lam{\Lambda}

\def\qed { Q.E.D. }
\def\6{\partial}

\def\dfq{\frac{\partial ^2 f}{\partial\bar q_i \partial  q_j}}

\swapnumbers
\newtheorem{theorem}{Theorem}[section]
\newtheorem{corollary}[theorem]{Corollary}
\newtheorem{lemma}[theorem]{Lemma}
\newtheorem{proposition}[theorem]{Proposition}
\newtheorem{claim}[theorem]{Claim}
\theoremstyle{definition}
\newtheorem{example}[theorem]{Example}
\newtheorem{definition}[theorem]{Definition}
\newtheorem{remark}[theorem]{Remark}
\newtheorem{proposition-definition}[theorem]{Proposition-Definition}

  \def\ci{{\cal I}}

 \makeatletter

\newcommand{\ps@verbit}{%
  \renewcommand{\@oddhead}{%
          \scriptsize
          {Plurisubharmonic functions on hypercomplex manifolds}
          \hfil\tiny {S. Alesker, M. Verbitsky}}
  \renewcommand{\@evenhead}{\@oddhead}
  \renewcommand{\@oddfoot}{\hfil\thepage\hfil}
  \renewcommand{\@evenfoot}{\@oddfoot}}

\@addtoreset{equation}{section} \@addtoreset{footnote}{section}

 \makeatother

\def\sh{S_\HH}
\begin{document}

\maketitle
\begin{abstract}

A hypercomplex manifold is a manifold equipped with a triple of
complex structures $I, J, K$ satisfying the quaternionic
relations. We define a quaternionic analogue of plurisubharmonic
functions on hypercomplex manifolds, and interpret these functions
geometrically as potentials of HKT (hyperk\"ahler with torsion)
metrics. We prove a quaternionic analogue of A.D. Aleksandrov and
Chern-Levine-Nirenberg theorems.
\end{abstract}

\tableofcontents

\section{Introduction.} The goal of this article is to introduce a
class of (continuous) quaternionic plurisubharmonic functions on
hypercomplex manifolds. We prove a version of A.D. Aleksandrov and
Chern-Levine-Nirenberg theorems for it. Then we present a
geometric characterization of smooth quaternionic strictly
plurisubharmonic functions as (local) potentials of so called
HKT-metrics on hypercomplex manifolds (HKT is the abbreviation of
HyperK\"ahler with Torsion). This interpretation is analogous to
the well known interpretation of smooth complex strictly
plurisubharmonic functions on complex manifolds as (local)
potentials of K\"ahler metrics.

The class of quaternionic plurisubharmonic functions on the flat
space $\HH^n$ was introduced by one of the authors in
\cite{alesker-bsm} and independently by G. Henkin \cite{henkin}
(unpublished). This class was studied in \cite{alesker-bsm},
\cite{alesker-ma}, \cite{alesker-valq}. Applications to the theory
of valuations on convex sets were obtained in \cite{alesker-valq}.
In this article we extend some of those definitions and results to
hypercomplex manifolds.

Other results related to quaternionic pluripotential theory on
hypercomplex manifolds were obtained by one of the authors
\cite{_V:reflexive_} (e.g. a quaternionic version of Sibony's
lemma \cite{_Sibony_} on extensions of positive currents and a
version of the Skoda-El Mir theorem).

Let us discuss the main results of this article in greater detail.

\begin{definition}
A {\itshape hypercomplex} manifold is a smooth manifold $X$
together with a triple $(I,J,K)$ of complex structures satisfying
the usual quaternionic relations:
$$IJ=-JI=K.$$
\end{definition}
\begin{remark}
(1) We will suppose in this article (in the opposite to much of
the literature on the subject) that the complex structures $I,J,K$
act on the {\itshape right} on the tangent bundle $TX$ of $X$.
This action extends uniquely to the right action of the algebra
$\HH$ of quaternions on $TX$.

(2) It follows that the dimension of a hypercomplex manifold $X$
is divisible by 4.

(3) Hypercomplex manifolds were introduced explicitly by Boyer
\cite{boyer}
\end{remark}

Let $(X^{4n},I,J,K)$ be a hypercomplex manifold. Let us denote by
$\Lam^{p,q}_I(X)$ the vector bundle of $(p,q)$-forms on the
complex manifold $(X,I)$. By the abuse of notation we will also
denote by the same symbol $\Lam^{p,q}_{I}(X)$ the space of
$C^\infty$-sections of this bundle.

Let
\begin{eqnarray}\label{l1}
\6\colon \Lambda_I^{p,q}(X)\to \Lam_I^{p+1,q}(X)
\end{eqnarray}
be the usual $\6$-differential of differential forms on the
complex manifold $(X,I)$.

Set
\begin{eqnarray}\label{l2}
\6_J:=J^{-1}\circ \bar \6 \circ J.
\end{eqnarray}

\begin{claim}[\cite{verbitsky-hkt}]\label{l3}
(1)$ J\colon\Lambda_I^{p,q}(X)\to\Lam_I^{q,p}(X).$

(2) $\6_J\colon \Lambda_I^{p,q}(X)\to\Lam_I^{p+1,q}(X).$

(3) $\6\6_J=-\6_J\6$.
\end{claim}

\begin{remark}
Claim \ref{l3} (1) is clear because $I$ and $J$ anticommute, and
Claim \ref{l3} (2) is directly implied by Claim \ref{l3} (1).
\end{remark}

\begin{definition}[\cite{verbitsky-hkt}]\label{l4}
Let $k=0,1,\dots,n$. A form $\ome\in \Lam^{2k,0}_I(X)$ is called
 \itshape{real} if
$$\overline{J\circ \ome}=\ome.$$
\end{definition}

We will denote the subspace of real $C^\infty$-smooth
$(2k,0)$-forms on $(X,I)$ by $\Lam^{2k,0}_{I,\RR}(X)$.
\begin{lemma}\label{l5}
Let $X$ be a hypercomplex manifold. Let $f\colon X\to \RR$ be a
smooth function. Then $\6\6_J f\in \Lam^{2,0}_{I,\RR}(X)$.
\end{lemma}
This lemma is proved in Section \ref{operators} as Lemma \ref{n5}.

\begin{definition}
Let $\ome\in \Lam^{2,0}_{I,\RR}(X)$. Let us say that $\ome$ is
non-negative (notation: $\ome\geq 0$) if
$$\ome(Y,Y\circ J)\geq 0$$
for any (real) vector field $Y$ on the manifold $X$. Equivalently,
$\omega$ is non-negative if $\omega(Z, \bar Z \circ J)\geq 0$ for
any $(1,0)$-vector field $Z$.
\end{definition}

\begin{definition}\label{l6}
A continuous function $$h:X\to \RR$$ is called quaternionic
plurisubharmonic if $\6\6_J h$ is a non-negative (generalized)
section of $\Lambda^{2,0}_{I,\RR}(X)$.
\end{definition}
\begin{remark}
The non-negativity in the generalized sense is discussed in detail
in Section \ref{vecbuncone}.
\end{remark}

Let us denote by $P'(X)$ the class of continuous quaternionic
plurisubharmonic functions on $X$. Let us denote by $P''(X)$ the
subclass of functions from $P'(X)$ with the following additional
property: a function $h\in P'(X)$ belongs to $P''(X)$ if and only
if any $x\in X$ has a neighborhood $U\ni x$ and a sequence
$\{h_N\}\subset P'(U)\cap C^2(U)$ such that
$h_N\overset{C^0(U)}{\to}h$ (where the convergence is understood
in sense of the uniform convergence on compact subsets of $U$).
Thus $P''(X)\subset P'(X)$.

We conjecture that $P'(X)=P''(X)$. This conjecture is true when
$X$ is an open subset of $\HH^n$.

The first main result of the article is the following theorem.
\begin{theorem}\label{l8}
Let $X$ be a hypercomplex manifold of (real) dimension $4n$. Let
$0<k\leq n$. For any $h_1,\dots,h_k\in P''(X)$ one can define a
non-negative generalized section of $\Lambda^{2k}_{I,\RR}$ denoted
by $\6\6_J h_1\wedge \dots\wedge\6\6_J h_k$ which is uniquely
characterized by the following two properties:

(1) if $h_1,\dots,h_k\in C^2(X)$ then the definition is clear;

(2) if $\{h_i^N\}\subset P''(X)$, $h_i^N\overset{C^0}{\to}h_i$ as
$N\to \infty$, $i=1,\dots, k$, then $h_i\in P''(X)$ and
$$\6\6_J h_1^N\wedge \dots\wedge \6\6_J h_k^N\to \6\6_J h_1\wedge \dots\wedge \6\6_J h_k$$
in the weak topology on measures.
\end{theorem}
This theorem is proved in Section \ref{qpsh} as Theorem \ref{2-9}.
\begin{remark}
Theorem \ref{l8} is a quaternionic analogue of a (real) result of
A.D. Aleksandrov \cite{aleksandrov2} and a (complex) result of
Chern-Levine-Nirenberg \cite{chern-levine-nirenberg}. When the
hypercomplex manifold $X$ is an open subset of the flat space
$\HH^n$ Theorem \ref{l8} was proved by one of the authors in
\cite{alesker-valq}, and in a special case of $k=n$ in
\cite{alesker-bsm}. Some applications of this theorem in the flat
case to the theory of valuations on convex sets were obtained in
\cite{alesker-valq}.
\end{remark}

In order to formulate the second main result we have to remind the
definition of an HKT-metric on a hypercomplex manifold $X$. Let
$g$ be a Riemannian metric on $X$. The metric $g$ is called
{\itshape quaternionic Hermitian} (or hyperhermitian) if $g$ is
invariant with respect to the group $SU(2)\subset \HH^*$ of
unitary quaternions.

Given a quaternionic Hermitian metric $g$ on a hypercomplex
manifold $X$, consider the differential form
$$\Ome:=\ome_J-\sqrt{-1}\ome_K$$
where $\ome_L(A,B):=g(A,B\circ L)$ for any $L\in \HH$ with
$L^2=-1$ and any vector fields $A,B$ on $X$. It is easy to see
that $\Ome$ is a $(2,0)$-form with respect to the complex
structure $I$.

\begin{definition}\label{def-hkt-metr}
The metric $g$ on $X$ is called HKT-metric if
$$\6 \Ome =0.$$
\end{definition}
\begin{remark}
HKT-metric on hypercomplex manifolds first were introduced by Howe
and Papadopoulos \cite{howe-papa}. Their original definition was
different but equivalent to Definition \ref{def-hkt-metr} (see
\cite{_Gra_Poon_}).
\end{remark}

Let us denote by $S_\HH(X)$ the vector bundle over $X$ such that
its fiber over a point $x\in X$ is equal to the space of
hyperhermitian forms on the tangent space $T_xX$ (see Definition
\ref{hyperform} in Section \ref{linalg}). Consider the map of
vector bundles
$$t\colon \Lam^{2,0}_{I,\RR}(X)\to S_\HH(X)$$
defined by $t(\eta)(A,A)=\eta(A,A\circ J)$ for any vector field
$A$ on $X$. The $t$ is an isomorphism of vector bundles (this was
proved in \cite{verbitsky-hkt}; see also Lemma \ref{n8} below).

The second main result is the following observation which provides
a geometric interpretation of the notion of quaternionic
(strictly) plurisubharmonic function on a hypercomplex manifold.
\begin{proposition}\label{m}
(1) Let $f$ be an infinitely smooth strictly plurisubharmonic
function on a hypercomplex manifold $(X,I,J,K)$. Then $t(\6\6_J
f)$ is an HKT-metric.

(2) Conversely assume that $g$ is an HKT-metric. Then any point
$x\in X$ has a neighborhood $U$ and an infinitely smooth strictly
plurisubharmonic function $f$ on $U$ such that $g=t(\6\6_J f)$ in
$U$.
\end{proposition}


This result is a direct consequence of a quaternionic version of
the local $\6\bar\6$-lemma well known for the complex manifolds.
We call this generalization local $\6\6_J$-lemma. It says as
follows.
\begin{proposition}\label{ddbar1}
Let $\Ome\in C^\infty(X,\Lam^{2,0}_{I,\RR})$. Then locally on $X$
the form $\Ome$ can be presented in a form
$$\Ome=\6\6_Jf$$
with $f$ being a $C^\infty$-smooth real valued function if and
only if $\6\Ome=0$.
\end{proposition}
The proof of this result is a rather straightforward application
of the main theorem of \cite{banos-swann} (which uses in turn
\cite{mamone-salamon}).


The article is organized as follows. In Section \ref{linalg} we
discuss some auxiliary constructions from quaternionic linear
algebra.

In Section \ref{operators} we discuss differential operators $\6$
and $\6_J$ on differential forms on general hypercomplex maifolds
and the so called Dirac operators on $\HH^n$.

In Section \ref{comparison} we make a comparison between
differential operators on the flat space $\HH^n$ and on general
hypercomplex manifolds; the goal is to rewrite some expressions on
$\HH^n$ in a language working in the more general setting of
hypercomplex manifolds.

In Section \ref{vecbuncone} we introduce a general notion of a
vector bundle with a cone in order to have a notion of positive
(with respect to this cone) section of the vector bundle.

In Section \ref{rempsh} we remind the definition and some results
on quaternionic plurisubharmonic functions on the flat space
$\HH^n$ following \cite{alesker-bsm}, \cite{alesker-valq}.

Section \ref{qpsh} contains the main definitions of this article
and the proof of the first main result Theorem \ref{l8}.

Section \ref{dhkt} describes the relation between quaternionic
plurisubharmonic functions and HKT-geometry. Namely we prove
Proposition \ref{m} and the local $\6\6_J$-lemma (Proposition
\ref{ddbar}).

{\bf Acknowledgements.} The first named author is grateful to G.
Henkin and V. Shevchishin for useful discussions. We thank the
referee for important remarks on the first version of the article.

\section{Some linear algebra.}\label{linalg} In this section we describe some
facts from linear algebra.

\begin{definition}
(1) Let $A=(a_{ij})$ be an $(n\times n)$-matrix with quaternionic
entries. Then $A$ is called {\itshape hyperhermitian} if
$$a_{ij}=\bar a_{ji}$$
where $\bar q$ is the usual conjugation of a quaternion $q$.

(2) A hyperhermitian matrix $A$ is called {\itshape non-negative
definite} (resp. {\itshape positive definite}) if for any $\xi\in
\HH^n\backslash\{0\}$ one has $\xi^*A\xi\geq 0$ (resp. $
\xi^*A\xi> 0$).
\end{definition}
Let $V$ be a right vector space over quaternions.
\begin{definition}\label{hyperform}
  A {\itshape hyperhermitian semilinear form}
on $V$ is a map $ a:V \times V \to \HH$ satisfying the following
properties:

(a) $a$ is additive with respect to each argument;

(b) $a(x,y \cdot q)= a(x, y) \cdot q$ for any $x,y \in V$ and any
$q\in \HH$;

(c) $a(x,y)= \overline{a(y,x)}$.
\end{definition}
\begin{remark}
It is easy to see that any hyperhermitian form $a$ on $\HH^n$ can
be written in the form $a(X,Y)=\sum_{i,j=1}^n\bar x_ia_{ij}y_j$
where $(a_{ij})$ is a uniquely determined hyperhermitian $n\times
n$- matrix.
\end{remark}
The space of hyperhermitian forms on $V$ we will denote by
$S_\HH(V)$.

For a quaternionic $n\times n$-matrix $A\in M_n(\HH)$ let us
denote by ${}^{\textbf{R}} A$ the {\itshape realization matrix} of
$A$ which is a real $4n\times 4n$-matrix. (To construct it,
consider $A$ as a matrix of a quaternionic transformation
$\HH^n\to \HH^n$. Identify $\HH^n\tilde \to \RR^{4n}$ in the
standard way. Then ${}^{\textbf{R}} A$ is the matrix of this
transformation with respect to the standard basis of $\RR^{4n}$.)

The following result is classical (see \cite{aslaksen} for the
references).
\begin{theorem}
There exists a polynomial $P$ defined on the space of all
hyperhermitian $n \times n$-matrices such that for any
hyperhermitian $n \times n$-matrix $A$ one has
$det({}^{\textbf{R}} A)= P^4(A)$ and $P(Id)=1$. $P$ is defined
uniquely by these two properties. Furthermore $P$ is homogeneous
of degree $n$ and has integer coefficients.
\end{theorem}
Thus for any hyperhermitian matrix $A$ the value $P(A)$ is a real
number, and it is called the {\itshape Moore determinant} of the
matrix $A$. The explicit formula for the Moore determinant  was
given by Moore \cite{moore} (see also the survey \cite{aslaksen}
and \cite{gelfand-retakh-wilson}). From now on the Moore
determinant of a matrix $A$ will be denoted by $det A$. This
notation should not cause any confusion with the usual determinant
of real or complex matrices due to part (i) of the next theorem.

\begin{theorem}

(i) The Moore determinant of any complex hermitian matrix
considered as quaternionic hyperhermitian matrix is equal to its
usual determinant.

(ii) For any hyperhermitian matrix $A$ and any quaternionic matrix
$C$
$$det (C^*AC)= detA \cdot det(C^*C).$$
\end{theorem}
For the proof we refer to \cite{alesker-bsm} though this result
was known earlier and is probably a folklore.
\begin{example}

(a) Let $A =diag(\lam_1, \dots, \lam _n)$ be a diagonal matrix
with real $\lam _i$'s. Then $A$ is hyperhermitian and the Moore
determinant $detA= \prod _i \lam_i$.

(b)  A general hyperhermitian $2 \times 2$ matrix $A$ has the form
 $$ A=  \left[ \begin {array}{cc}
                     a&q\\
                \bar q&b\\
                \end{array} \right] $$
where $a,b \in \RR, \, q \in \HH$. Then $det A =ab - q \bar q
(=ab- \bar q q)$.
\end{example}

Let $V$ be a right $\HH$-module of quaternionic dimension $n$. Let
$0\leq k\leq n$ be an integer. Let us denote by
$$\Lam^{2k,0}_I(V):=\wedge^{2k}_{\CC}(V_I^*)$$
where $V_I$ is the space $V$ equipped with the complex structure
$I$, $V_I^*$ is its dual. Note that
$$J\colon V_I\to V_I$$
is an anti-linear map (namely $J(x\cdot \lam)=J(x)\cdot
\bar\lam$). It induces an anti-linear {\itshape involution}
$$J\colon \Lam^{2k,0}_I(V)\to \Lam^{2k,0}_I(V).$$
\def\lir{\Lam^{2k,0}_{I,\RR}(V)}
Let us denote by $\lir$ the real subspace fixed by this
involution.
\begin{lemma}\label{g1}
The natural representation of the group $GL_n(\HH)$ in $\lir$ is
absolutely irreducible (in particular it is irreducible). The
complexification of this representation has highest weight
$(\underbrace{0,\dots,0}_{2(n-k) \mbox{
times}},\underbrace{-1,\dots,-1}_{2k \mbox{ times}})$ as a
representation of $GL_n(\HH)\otimes_\RR \CC=GL_{2n}(\CC)$.
\end{lemma}
To prove Lemma \ref{g1} we will need the following elementary
lemma which is in fact a special case of Hilbert 90 Theorem.
\begin{lemma}\label{g2}
Let $W$ be a complex vector space. Let $\sigma\colon W\to W$ be an
anti-linear involution of $W$. Then
$$W=W^\sigma\oplus W^\sigma\cdot \sqrt{-1}$$
where $W^\sigma$ is the (real) subspace of $\sigma$-fixed vectors.
\end{lemma}
{\bf Proof} of Lemma \ref{g1} assuming Lemma \ref{g2}. By Lemma
\ref{g2} we have
$$\Lam^{2k,0}_I(V)=\lir\oplus \lir\cdot I.$$ Hence
$\Lam^{2k,0}_I(V)$ is the complexification of the representation
of $GL_n(\HH)$ in $\lir$. The complexification of the group
$GL_{n}(\HH)$ is the group $GL_{2n}(\CC)$. But the representation
of the group $GL_{2n}(\CC)$ in
$\Lam^{2k,0}_I(V)=\wedge^{2k}_\CC(V_I^*)$ is irreducible with
highest weight $(\underbrace{0,\dots,0}_{2(n-k) \mbox{
times}},\underbrace{-1,\dots,-1}_{2k \mbox{ times}})$ (this is a
basic fact from representation theory). \qed

{\bf Proof} of Lemma \ref{g1}. Set
$$W^{\sigma'}:=\{x\in W|\, \sigma x=-x\}.$$
Then clearly $W=W^\sigma\oplus W^{\sigma'}$. Obviously $x\in
W^\sigma$ if and only if $x\cdot I\in W^{\sigma'}$, hence
$W^{\sigma'}=W^\sigma\cdot I$. \qed

\begin{lemma}\label{n6}
Let $\ome\in \Lam^{2,0}_{I,\RR}(V)$. Define a quadratic form
$$B(X,X)=\ome(X,X\circ J).$$
Then $B$ is hyperhermitian, i.e. $B\in \sh(V)$.
\end{lemma}
{\bf Proof.} First let us check that $B$ is real valued. Indeed
\begin{eqnarray*}
\overline{B(X,X)}=\overline{\ome(X,X\circ J)}=\ome(X\circ J,X\circ
J^2)=\\
\ome(X,X\circ J)=B(X,X).
\end{eqnarray*}

Let us check that $B$ is invariant under $I,J,K$:
\begin{eqnarray*}
B(X\circ I,X\circ I)=\ome(X\circ I,X\circ IJ)=\\
-\ome(X\circ I,(X\circ J)\circ I)=\ome(X,X\circ J)=B(X,X);\\
B(X\circ J,X\circ J)=\ome(X\circ J,-X)=\ome(X,X\circ J)=B(X,X);\\
B(X\circ K,X\circ K)=B(X\circ IJ,X\circ IJ)=B(X\circ I,X\circ
I)=B(X,X).
\end{eqnarray*}
\qed

Thus Lemma \ref{n6} defines a map
\begin{eqnarray}\label{n7}
t\colon \Lam_{I,\RR}^{2,0}(V)\to \sh(V).
\end{eqnarray}
\begin{lemma}[\cite{verbitsky-hkt}]\label{n8}
The map $t$ defined in (\ref{n7}) is an isomorphism. The inverse
map is given by
\begin{eqnarray}\label{n8.5}
(t^{-1}g)(X,Y)=-(g(X,Y\circ J)-\sqrt{-1}g(X,Y\circ
K)).\end{eqnarray}
\end{lemma}
{\bf Proof.} Let us consider the map
$$\phi\colon S_\HH(V)\to \Lam^{2,0}_{I,\RR}(V)$$
defined by $(\phi (g))(X,Y)=-(g(X,Y\circ J)-\sqrt{-1}g(X,Y\circ
K))$. We will show that $\phi$ is the inverse of $t$. But first
let us check that indeed for any $g\in S_{\HH}(V)$ we have
\begin{eqnarray}\label{n8.6}
\phi (g)\in \Lam^{2,0}_{I,\RR}(V).
\end{eqnarray}
To prove that $\phi(g)\in \Lam^{2,0}_I(V)$ it is enough to check
that
\begin{eqnarray*}
\frac{d}{d\theta}|_0(\phi(g))(X\circ e^{\theta I},Y\circ e^{\theta
I})=2\sqrt{-1}(\phi(g))(X,Y).
\end{eqnarray*}
We have
\begin{multline*}
\frac{d}{d\theta}|_0(\phi(g))(X\circ e^{\theta I},Y\circ
e^{\theta I})=\phi(g)(X\circ I,Y)+\phi(g)(X,Y\circ I) \\
=-(g(X\circ I,Y\circ J)-\sqrt{-1}g(X\circ I,Y\circ K)+g(X,Y\circ
IJ)-\sqrt{-1}g(X,Y\circ IK)) \\
= -(g(X,Y\circ K)+\sqrt{-1}g(X,Y\circ J)+g(X,Y\circ
K)+\sqrt{-1}g(X,Y\circ J))\\
= 2\sqrt{-1}(\phi(g))(X,Y)
\end{multline*}
Let us show that $\phi(g)$ is real. We have
\begin{multline*}
\overline{(\phi(g))(X\circ J,Y\circ J)}\\=-(g(X\circ J,Y\circ
J^2)+\sqrt{-1} g(X\circ J,Y\circ JK))=(\phi(g))(X,Y).
\end{multline*}
Thus (\ref{n8.6}) is proved.

Let $g\in S_\HH(V)$. Then
\begin{eqnarray*}
(t\circ \phi)(g)(X,X)=-(g(X,X\circ J^2)-\sqrt{-1}g(X,X\circ
JK))=g(X,X).
\end{eqnarray*}
Hence
\begin{eqnarray}\label{n8.7}
t\circ \phi =id.
\end{eqnarray}
Let us check that $\phi\circ t=id$. Let $\eta\in
\Lam^{2,0}_{I,\RR}(V)$. Set $g:=t(\eta)$. Then
\begin{eqnarray*}
g(X,Y)=\frac{1}{2}(g(X+Y,X+Y)-g(X,X)-g(Y,Y))=\\
\frac{1}{2}(\eta(X+Y,(X+Y)\circ J)-\eta(X,X\circ J)-\eta(Y,Y\circ
J)).
\end{eqnarray*}
Then
\begin{multline*}
(\phi\circ t)(\eta)(X,Y)=-(g(X,Y\circ J)-\sqrt{-1}g(X,Y\circ K))\\
= -\frac{1}{2} ( g(X+Y\circ J,X+Y\circ
J)-g(X,X)-g(Y,Y)-\\\sqrt{-1}(g(X+Y\circ K,X+Y\circ
K)-g(X,X)-g(Y,Y)))\\ = -\frac{1}{2} ( \eta(X+Y\circ J,X\circ
J-Y)-\eta(X,X\circ J)-\eta(Y,Y\circ J)-\\\sqrt{-1}(\eta(X+Y\circ
K,X\circ J-Y\circ I)-\eta(X,X\circ J)-\eta(Y,Y\circ J))).
\end{multline*}
Opening by bilinearity and making cancellations the last
expression becomes
\begin{multline*}
-\frac{1}{2} (-\eta(X,Y)-\eta(X\circ J,Y\circ J)- \\
\sqrt{-1}(-\eta(X,Y\circ I)+\eta(Y\circ K,X\circ J)-\eta(Y\circ
IJ,Y\circ I)+\eta(Y\circ K,Y\circ I)))\\= \frac{1}{2}
(\eta(X,Y)+\overline{\eta(X,Y)}-\sqrt{-1}(\eta(X,Y\circ
I)-\overline{\eta(Y\circ I,X)})\\= Re\,\eta(X,Y)-\sqrt{-1}
Re\,\eta(X,Y\circ I)=\eta(X,Y).
\end{multline*}
Thus Lemma \ref{n8} is proved. \qed

Now we are going to define convex cones of strongly and weakly
positive elements in $\lir$, $0\leq k\leq n$. The exposition is
analogous to the complex case as in Harvey \cite{harvey} (see also
Lelong \cite{lelong}). First observe that $\Lam^{2n,0}_{I,\RR}(V)$
is a real one-dimensional real space. It is canonically oriented.
Let us denote the (closed) half line of positive elements in
$\Lam^{2n,0}_{I,\RR}(V)$ by $\Lam^{2n,0}_{I,\RR}(V)_{\geq 0}$.
\begin{definition}\label{g1.5}
(1) An element $\eta\in\lir$ is called {\itshape strongly
positive} if it can be presented as a finite sum of elements of
the form $f^*\xi$ where $f\colon V\to U$ is a morphism of right
$\HH$-modules, $\dim_\HH U=k$, $\xi\in
\Lam^{2k,0}_{I,\RR}(U)_{\geq 0}$.

(2) An element $\eta\in \lir$ is called {\itshape weakly positive}
if for any strongly positive element $\zeta\in
\Lam^{2(n-k),0}_{I,\RR}(V)$ one has $\eta\wedge \zeta\in
\Lam^{2n,0}_{I,\RR}(V)_{\geq 0}$.
\end{definition}
It is clear that strongly and weakly positive elements form convex
cones. Let us denote by $C^k(V)$ (resp. $K^k(V)$) the cone of
strongly (resp. weakly) positive elements.
\begin{remark}
(1) Clearly we have $(\lir)^*=\Lam^{2(n-k),0}_{I,\RR}(V)\otimes
\Lam^{2n,0}_{I,\RR}$. Then the closure $\overline{C^k(V)}$ is the
cone dual to the cone $K^{n-k}\otimes\Lam^{2n,0}_{I,\RR}(V)_{\geq
0}$. (The duality of cones is understood in the standard sense:
for a convex cone $K$ in a vector space $W$ one defines the dual
cone $K^\circ:=\{y\in W^*|\, y(x)\geq 0 \,\,\forall x\in K\}$.)

(2)We will see below in Propositions \ref{g3}(2), \ref{g7} that
$\overline{C^2(V)}=K^2(V)$ and this cone coincides with the cone
$\{\eta\in \Lam^{2,0}_{I,\RR}(V)|\, \eta(A,A\circ J)\geq 0
\,\,\forall A\in V\}$.
\end{remark}
Let us state some basic properties of the cones $C^k(V), K^k(V)$.
\begin{proposition}\label{g3}
(1) $C^k(V)\subset K^k(V)$.

(2) $C^k(V)\wedge C^l(V)\subset C^{k+l}(V)$.

(3) For $k=0,1,n-1,n$ $$C^k(V)=K^k(V).$$

(4) The cones $C^k(V)$ and $K^k(V)$ have non-empty interior.
\end{proposition}
This proposition was proved in \cite{alesker-valq} in a somewhat
different language. Now we will do this comparison of languages.

First let us describe the relevant linear algebraic constructions
from \cite{alesker-valq}, Section 2. Remind that we denote by
$S_\HH(V)$ the space of hyperhermitian forms on $V$. Recall
$\dim_\HH V=n$.
\begin{proposition}\label{omega-def}[\cite{alesker-valq}, Section 2]
Let $0\leq k\leq n$ be an integer. The space $Sym^k_\RR(S_\HH(V))$
has a unique quotient denoted by $\Ome^{k,k}(V)$ such that the
complexification of the natural representation of $GL_n(\HH)$ in
$\Ome^{k,k}(V)$ is irreducible and has highest weight
$(\underbrace{0,\dots,0}_{2(n-k)\mbox{
times}},\underbrace{-1,\dots,-1}_{2k\mbox{ times}})$ as a
representation of $GL_{2n}(\CC)=GL_n(\HH)\otimes _\RR \CC$.
\end{proposition}
Let $$p_k\colon Sym^k(S_\HH(V))\to \Ome^{k,k}(V)$$ be the
canonical projection.

Define $$\Ome^\bullet(V):=\oplus_{k=0}^n\Ome^{k,k}(V).$$ Let us
describe the algebra structure on $\Ome^\bullet(V)$ following
\cite{alesker-valq}. Consider the composition $\mu$ of maps
\begin{eqnarray*}
Sym^k(S_\HH(V)\otimes Sym^l(S_\HH(V))\to
Sym^{k+l}(S_\HH(V))\overset{p_{k+l}}{\to}\Ome^{k+l,k+l}(V).
\end{eqnarray*}
It was shown in Proposition 2.1.11 of \cite{alesker-valq} that
$\mu$ factorizes (uniquely) via $p_k\otimes p_l$, namely there
exists a unique map
$$m\colon \Ome^{k,k}(V)\otimes \Ome^{l,l}(V)\to \Ome^{k+l,k+l}(V)$$
which makes the following diagram commutative:
\[\begin{CD}
Sym^k(S_\HH(V))\otimes Sym^l(S_\HH(V)) @>>>
Sym^{k+l}(S_\HH(V)) \\
@VVV @VVV\\
\Ome^{k,k}(V)\otimes\Ome^{l,l}(V) @>>m> \Ome^{k+l,k+l}(V)
\end{CD}
\]
This map $m$ defines the product on $\Ome^\bullet (V)$.
\begin{proposition}[\cite{alesker-valq}, Theorem 2.1.13]\label{g4}
The correspondence \[ V\mapsto \Ome^\bullet(V)\] is a
contravariant functor from the category of finite dimensional
$\HH$-modules to the category of finite dimensional commutative
associative graded algebras. For a fixed $V$ the graded algebra
$\Ome^\bullet(V)$ satisfies the Poincar\'e duality.
\end{proposition}
Recall that we have the isomorphism
$$t\colon \Lam^{2,0}_{I,\RR}\tilde\to S_\HH(V).$$
Clearly $t$ commutes with the natural action of $GL_n(\HH)$. Fix
an integer $0\leq k\leq n$. We have the canonical map
$$\gamma_k\colon Sym^k(\Lam^{2,0}_{I,\RR}(V))\to
\Lam^{2k,0}_{I,\RR}(V)$$ given by $\eta_1\otimes \dots\otimes
\eta_k\mapsto \eta_1\wedge\dots\wedge\eta_k$.
\begin{lemma}\label{g5}
There exists a unique map
$$\tau_k\colon \Lam^{2k,0}_{I,\RR}(V)\to \Ome^{k,k}(V)$$
which makes the following diagram commutative:
\[ \begin{CD}
Sym^k(\Lam^{2,0}_{I,\RR}(V)) @>{Sym^k t}>>Sym^{k}(S_\HH(V))
 \\
@V{\gamma_k}VV @VV{p_k}V \\
\Lam^{2k,0}_{I,\RR}(V) @>>{\tau_k}> \Ome^{k,k}(V).
\end{CD}
\]
the map $\tau_k$ is an isomorphism and commutes with the action of
$GL_n(\HH)$.
\end{lemma}
{\bf Proof.} This lemma follows immediately from the following
fact:

(a) $t$ is isomorphism of $GL_n(\HH)$-modules;

(b) $Sym^k(S_\HH(V))\otimes_\RR \CC$ is multiplicity free as
$GL_{2n}(\CC)=GL_n(\HH)\otimes_\RR\CC$-module (it is Proposition
2.1.7 of \cite{alesker-valq} combined with Lemma 2.1.4 of
\cite{alesker-valq}; Proposition 2.1.7 of \cite{alesker-valq} is
due to \cite{howe}, Proposition 2);

(c) $\lir\otimes_\RR\CC$ and $\Ome^{k,k}(V)\otimes_\RR\CC$ are
irreducible $GL_{2n}(\CC)$-modules of the same highest weight.

\qed

{}From the construction of the product on $\Ome^\bullet(V)$ and
the isomorphisms $\tau_k\colon \lir\tilde\to \Ome^{k,k}(V)$ one
easily has the following result.
\begin{proposition}\label{g6}
\[ \oplus_{k=0}^n \tau_k\colon
\oplus_{k=0}^n\Lam^{2k,0}_{I,\RR}(V)\tilde\to \Ome^\bullet(V)\] is
an isomorphism of graded algebras where the algebra structure on
\[ \oplus_{k=0}^n\Lam^{2k,0}_{I,\RR}(V)\] is the usual wedge product.
\end{proposition}

In \cite{alesker-valq}, Subsection 2.2, we have defined the cones
$ C^k(V)$ (resp. $K^k(V)$) in $\Ome^{k,k}(V)$ of strongly (resp.
weakly) positive elements. Definition \ref{g1.5} of cones in
$\lir$ in this article is obtained by applying the isomorphism
$\tau_k$ to those cones in $\Ome^{k,k}(V)$. Hence Proposition
\ref{g3} follows immediately from the corresponding properties of
the cones in $\Ome^\bullet(V)$ (see Propositions 2.2.2-2.2.5 in
\cite{alesker-valq}).
\begin{proposition}\label{g7}
The cone $C^2(V)\subset \Lam^{2,0}_{I,\RR}(V)$ is equal to the set
$\{\eta\in \Lam^{2,0}_{I,\RR}(V)|\, \eta(A,A\circ J)\geq 0
\,\,\forall A\in V\}$.
\end{proposition}
{\bf Proof.} Consider the cone $t(C^2(V))\subset S_\HH(V)$. This
is the cone of strongly positive elements in $S_\HH(V)$ considered
in \cite{alesker-valq}. In the proof of Proposition 2.2.4 of
\cite{alesker-valq} it was shown that this cone coincides with the
cone of non-negative definite hyperhermitian matrices. But
$t(\eta)(A,A)=\eta(A,A\circ J)$. Proposition \ref{g7} is proved.
\qed

\section{Differential operators  on hypercomplex\\
manifolds.}\label{operators} Let $(X^{4n},I,J,K)$ be a
hypercomplex manifold. Remind that we denote by $\Lam^{p,q}_I$ the
vector bundle of $(p,q)$-forms on the complex manifold $(X,I)$. By
the abuse of notation we will denote by this symbol also the space
of $C^\infty$-sections of this bundle.

Let
\begin{eqnarray}\label{n1}
\6\colon \Lambda_I^{p,q}(X)\to \Lam_I^{p+1,q}(X)
\end{eqnarray}
be the usual $\6$-differential of differential forms on the
complex manifold $(X,I)$. Set
\begin{eqnarray}\label{n2}
\6_J:=J^{-1}\circ \bar \6 \circ J.
\end{eqnarray}
\begin{claim}[\cite{verbitsky-hkt}; see also Claim \ref{l2}]\label{n3}
(1)$ J\colon\Lambda_I^{p,q}(X)\to\Lam_I^{q,p}(X).$

(2) $\6_J\colon \Lambda_I^{p,q}(X)\to\Lam_I^{p+1,q}(X).$

(3) $\6\6_J=-\6_J\6$.
\end{claim}
\begin{definition}[\cite{verbitsky-hkt}]\label{n4}
Let $k=0,\dots,n$. A form $\ome\in \Lam^{2k,0}_I(X)$ is called
 \itshape{real} if if it is real pointwise:
$$\overline{J\circ \ome}=\ome.$$
\end{definition}
\begin{lemma}\label{n5}
Let $f\colon X\to \RR$ be a real valued smooth function.  Then
$\6\6_J f\in \Lam_I^{2,0}(X)$ is real.
\end{lemma}
{\bf Proof.}
\begin{eqnarray*}
\overline{J\circ (\6\6_J f)} =\overline{J\circ \6 \circ
J^{-1}\circ \bar \6 f}=\\
J\circ \bar \6 \circ J^{-1}\circ \overline{\bar \6 f}=J\circ \bar
\6 \circ J^{-1}\circ \6\bar f=J\circ \bar \6 \circ
J^{-1}\circ \6 f=\\
-\6_J\6 f=\6\6_J f
\end{eqnarray*}
where the last equality is by Claim \ref{n3}(3). \qed

{\bf Notation:} (1) We will denote the subspace of real
$C^\infty$-smooth $(2k,0)$-forms on $(X,I)$ by
$\Lam^{2k,0}_{I,\RR}(X)$.

\def\sh{S_\HH}
(2) We will denote by $S_{\HH}(X)$ the vector bundle over $X$ with
fiber over a point $x\in X$ equal to the space of hyperhermitian
quadratic forms on the tangent space $T_xX$.


\def\db{\frac{\partial}{\partial \bar q}}
\def\dq{\frac{\partial}{\partial  q}}

On the flat space $\HH^n$ one can introduce so called Dirac
operators. Let us describe them. We will write a quaternion $q$ in
the usual form
$$q= t+ x\cdot i +y\cdot j+ z\cdot k $$
where $t,\, x,\, y,\, z$ are real numbers, and $i,\, j,\, k$
satisfy the usual relations
$$i^2=j^2=k^2=-1, \, ij=-ji=k,\, jk=-kj=i, \, ki=-ik=j.$$

The Dirac operator $\frac {\partial}{\partial \bar q}$ is defined
as follows. For any $\HH$-valued function $f$
\begin{eqnarray}\label{dqbar}
\db f:=\frac{\partial f}{\partial  t}  + i \frac{\partial
f}{\partial x} + j \frac{\partial f}{\partial y} + k
\frac{\partial f}{\partial  z}.
\end{eqnarray}

Let us also define the operator $\dq$:
\begin{eqnarray}\label{dq}
\dq f:=\overline{ \db \bar f}= \frac{\partial f}{\partial  t}  -
 \frac{\partial f}{\partial x}  i-
 \frac{\partial f}{\partial  y} j-
\frac{\partial f}{\partial  z}  k.
\end{eqnarray}

It is easy to see that on $\HH^n$
$$\left[\frac{\6}{\6 q_i},\frac{\6}{\6\bar q_j}\right]=0.$$
It is easy to see that if $f\colon \HH^n\to \RR$ is a $C^2$-smooth
function then the matrix $\left(\frac{\6^2f}{\6\bar q_i\6
q_j}\right)$ is hyperhermitian (see \cite{alesker-bsm}).

\section{Comparison with the flat case.}\label{comparison}
Consider $\HH^n$ as a right $\HH$-vector space equipped with the
standard coordinate system.

\begin{proposition}\label{n9}
Let $f\colon \HH^n\to \RR$ be a real valued smooth function. Then
\begin{eqnarray}\label{n10}
t(\6\6_J f)=\frac{1}{4} \left(\dfq\right).
\end{eqnarray}
\end{proposition}
It is enough to prove the formula (\ref{n10}) pointwise. Since all
the expressions involved are equivariant under translations it is
enough to prove (\ref{n10}) at 0. It is enough to show that the
hyperhermitian quadratic forms $(t(\6 \6_Jf))|_0$ and $\left(
\dfq\right)(0)$ coincide on each quaternionic line . Since all the
operators in (\ref{n10}) are equivariant under the group
$GL_n(\HH)$ it is enough to assume that this quaternionic line is
equal to $\{(q,0,\dots,0)|\, q\in \HH\}$. Thus we may assume that
$n=1$. Then we have
\begin{eqnarray*}
\frac{\6^2f}{\6\bar q\6  q}= \frac{\6^2f}{\6 t ^2}+\frac{\6^2f}{\6
x^2}+\frac{\6^2f}{\6 y^2}+ \frac{\6^2f}{\6 z^2}
\end{eqnarray*}
when $q=t+xI+yJ+zK$.

Let us compute the left hand side in (\ref{n10}). Let us identify
$\HH^1\simeq \CC^2$ as follows: for $q=t+xI+yJ+zK=(t+xI)+J(y-zI)$
define
\begin{eqnarray*}
z_1:=t+xI,\\
z_2:=y-zI.
\end{eqnarray*}
Thus $q=z_1+Jz_2$. We have
\begin{eqnarray}\label{n11}
\frac{\6}{\6\bar z_1}=\frac{1}{2} (\frac{\6}{\6
t}+\sqrt{-1}\frac{\6}{\6 x});\\
\frac{\6}{\6\bar z_2}=\frac{1}{2} (\frac{\6}{\6
y}-\sqrt{-1}\frac{\6}{\6 z}).
\end{eqnarray}

\begin{claim}\label{n12}
\begin{eqnarray*}
J^{-1}\circ d\bar z_1=dz_2,\,J^{-1}\circ d z_1=d\bar z_2;\\
J^{-1}\circ d\bar z_2=-dz_1,\, J^{-1}\circ d z_2=-d\bar z_1.
\end{eqnarray*}
\end{claim}
{\bf Proof} is a straightforward computation. \qed

 Next we have
\begin{align*}
-\6_J\6 f=&-\6_J\left(\frac{\6 f}{\6 z_1}dz_1+\frac{\6
  f}{\6 z_2}dz_2\right)= -J^{-1}\circ \bar\6\left(\frac{\6 f}{\6 z_1}(J\circ
dz_1)+\frac{\6 f}{\6 z_2}(J\circ dz_2)\right) \\
= & J^{-1}\circ \bar\6 \left(\frac{\6 f}{\6 z_1}d\bar
  z_2-\frac{\6 f}{\6 z_2}d\bar z_1\right)
\\=& J^{-1}\circ\left(\frac{\6^2 f}{\6 z_1\6\bar z_1} d\bar
  z_1\wedge d\bar z_2-\frac{\6^2 f}{\6 z_2\6\bar z_2}
  d\bar z_2\wedge d\bar z_1\right)\\= &
\left(\frac{\6^2 f}{\6 z_1\6\bar z_1}+\frac{\6^2 f}{\6 z_2\6\bar
z_2}\right)\cdot J^{-1}\circ (d\bar z_1\wedge d\bar z_2)\\ =
&\frac{1}{4} \left(\frac{\6^2f}{\6 t ^2}+\frac{\6^2f}{\6
x^2}+\frac{\6^2f}{\6 y^2}+ \frac{\6^2f}{\6 z^2}\right)(dz_1\wedge
dz_2) \\ =& \frac{1}{4}\frac{\6^2 f}{\6\bar q\6 q}(dz_1\wedge
dz_2).
\end{align*}

Hence $$t(\6\6_J f)=\frac{1}{4}\frac{\6^2 f}{\6\bar q\6 q}\cdot
t(dz_1\wedge dz_2).$$ Hence in order to finish the proof of
Proposition \ref{n9} it remains to prove the following claim.
\begin{claim}\label{n13}
$$t(dz_1\wedge dz_2)=dt^2+dx^2+dy^2+dz^2.$$
\end{claim}
{\bf Proof} is a straightforward computation. \qed

\begin{remark}\label{n14}
The form $dz_1\wedge dz_2\wedge\dots \wedge dz_{2n-1}\wedge
dz_{2n}$ belongs to $\Lam^{2n,0}_{I,\RR}$ on $\HH^n$.
\end{remark}
Let us denote by $\underline\RR$ the trivial real line bundle over
$\HH^n$. Let  us consider the isomorphism of line bundles
\begin{eqnarray*}
F\colon \Lam^{2n,0}_{I,\RR}\to \underline\RR
\end{eqnarray*}
defined by $$F(\eta):=\frac{1}{n!}\frac{\eta}{dz_1\wedge
dz_2\wedge\dots \wedge dz_{2n-1}\wedge dz_{2n}}.$$ Let us denote
by
\begin{eqnarray*}
\kappa\colon Sym^{n}(\Lam^{2,0}_{I,\RR})\to \Lam^{2n,0}_{I,\RR}
\end{eqnarray*}
the natural map of vector bundles given by
$$\kappa(\eta_1\otimes \dots\otimes
\eta_n):=\eta_1\wedge\dots\wedge\eta_n.$$

\def\smnt{sym^n t}
Let us denote by $\det(X_1,\dots,X_n)$ the mixed determinant of
hyperhermitian matrices $X_1,\dots,X_n\in S_\HH(\HH^n)$. By
definition, it is a polarization of the Moore determinant. More
precisely, the mixed determinant is a map $Sym^n(S_\HH(\HH^n))\to
\RR$ which is uniquely characterized by the following property:
for any hyperhermitian matrix $Y\in S_\HH(\HH^n)$,
$\det(Y,\dots,Y)$ is equal to the Moore determinant $\det Y$.

\begin{proposition}\label{n15}
The following diagram is commutative
\begin{equation}
\label{n16}
\begin{CD}
Sym^n(\Lam^{2,0}_{I,\RR}) @>{Sym^n t}>> Sym^n(S_\HH)\\
@V\kappa VV @VV{det} V\\
\Lam^{2n,0}_{I,\RR} @>>F> \underline{\RR}
\end{CD}
\end{equation}
\end{proposition}
{\bf Proof.} It is enough to check the commutativity of this
diagram fiberwise. On each fiber we have the natural action of the
group $SL_n(\HH)$ and all the maps commute with the action of this
group (note that the action of this group on the spaces of the
bottom line is trivial). By Lemma 2.1.4 of \cite{alesker-valq} and
Proposition 2 of \cite{howe} (which is cited in Proposition 2.1.7
of \cite{alesker-valq}) the complexified representation of
$GL_n(\HH)$ (and hence of $SL_n(\HH)$) in fibers of $Sym^n(S_\HH)$
is multiplicity free. Since the map $t$ commutes with the action
of $GL_n(\HH)$, the complexified representation of $SL_n(\HH)$ in
fibers of $\Lam^{2,0}_{I,\RR}$ is also multiplicity free. Thus the
maps $\kappa,\, det$ are characterized uniquely up to a constant
by the property that they commute with the action of $SL_n(\HH)$.
Hence the diagram (\ref{n16}) must be commutative up to a
constant. To check the constant let us take $(dz_1\wedge
dz_2)\otimes \dots \otimes (dz_{2n-1}\wedge dz_{2n})$. We have
\begin{eqnarray*}
det (sym^n t((dz_1\wedge dz_2)\otimes \dots \otimes
(dz_{2n-1}\wedge dz_{2n})))=\\
det \left(\left[ \begin{array}{cccc}
                                       1&0&\dots&0\\
                                       0&0&\dots&0\\
                                       \multicolumn{4}{c}{\dotfill}\\
                                       0&0&\dots&0
                                       \end{array}\right],\left[ \begin{array}{cccc}
                                       0&0&\dots&0\\
                                       0&1&\dots&0\\
                                       \multicolumn{4}{c}{\dotfill}\\
                                       0&0&\dots&0
                                       \end{array}\right],\dots,\left[ \begin{array}{cccc}
                                       0&0&\dots&0\\
                                       0&0&\dots&0\\
                                       \multicolumn{4}{c}{\dotfill}\\
                                       0&0&\dots&1
                                       \end{array}\right]
                                       \right)=\\
\frac{1}{n!}\frac{\6^n}{\6\lam_1\dots \6\lam_n}|_0
det\left[\begin{array}{cccc}
                                 \lam_1&0&\dots&0\\
                                 0&\lam_2&\dots&0\\
                                 \multicolumn{4}{c}{\dotfill}\\
                                  0&0&\dots&\lam_n
                                  \end{array}\right]=\frac{1}{n!}.
\end{eqnarray*}
On the other hand
\begin{multline*}
\kappa((dz_1\wedge dz_2)\otimes \dots \otimes (dz_{2n-1}\wedge
dz_{2n}))\\=(dz_1\wedge dz_2)\wedge \dots \wedge (dz_{2n-1}\wedge
dz_{2n}).
\end{multline*}
Proposition \ref{n15} is proved. \qed

From Propositions \ref{n9} and \ref{n15} we immediately get the
following corollary.
\begin{corollary}\label{n17}
$$(\6\6_J f)^n=\frac{n!}{4^n}det(\frac{\6^2f}{\6\bar q_i\6
q_j})dz_1\wedge dz_2\wedge\dots \wedge dz_{2n-1}\wedge
dz_{2n}.$$
\end{corollary}
\section{Vector bundles with a cone.}\label{vecbuncone}
\begin{definition}\label{1-1} Let $E\to X$ be a finite dimensional real
vector bundle over a manifold $X$. We say that $E$ has a {\itshape
cone} if at each fiber $E_x,\, x\in X$, we are given a convex cone
$C_x\subset E_x$ with the following property: any point $x\in X$
has a neighborhood $U$, a trivialization $\phi: E|_U\tilde \to
X\times V$, and a convex cone $C\subset V$ such that for any $x\in
U$ one has $\phi(C_x)=\{x\}\times C$.
\end{definition}

Note that if a vector bundle $E$ has a cone then the dual bundle
$E^*$ has a cone (which is dual to the cone of $E$).

\def\ox{|\omega_X|}
Let us denote by $|\omega_X|$ the real line bundle of densities on
$X$. It is canonically oriented. Hence if a bundle $E$ has a cone
then naturally $E\otimes \ox$ has a cone.

\begin{definition}\label{1-2}
(1) Let $E$ be a vector bundle with a cone. A continuous section
$\gamma \in C(X,E)$ is called {\itshape non-negative} if at any
point $x\in X$ $\gamma(x)\in \bar{C_x}$ (where $\bar C_x$ denotes
the closure of the cone $C_x$).

(2) A generalized section $\gamma \in C^{-\infty}(X,E)$ is called
{\itshape non-negative} if for any non-negative section $\phi\in
C^\infty_0(X,E^*\otimes \ox)$ one has $\gamma(\phi)\geq 0$.
\end{definition}

It is easy to see that a continuous section $\gamma\in C(X,E)$ is
non-negative if and only if it is non-negative as a generalized
section.

\begin{definition}\label{3}
Let $E$ be a vector bundle. An {\itshape $E$-valued measure} is a
continuous linear functional $C_0(X,E^*)\to \RR$.
\end{definition}
\begin{proposition}\label{1-3}
Let $E$ be a vector bundle with a cone such that the dual bundle
$E^*$ has a cone with non-empty interior at each point. Then any
non-negative generalized section of $E$ is an $E\otimes
|\ome_X|^*$-valued measure.
\end{proposition}
{\bf Proof.} This proposition is essentially well known. (When $E$
is the trivial line bundle this is proved in
\cite{generalized-functions}, Ch. II.) Using partition of unity
the proof reduces to the case of a trivial bundle $E$.
\qed

Let us describe the now the main examples of bundles with a cone
which will be used in this article. Let $(X^{4n},I,J,K)$ be a
hypercomplex manifold. Let $0\leq k\leq n$ be an integer.  The
bundle $\Lam^{2k,0}_{I,\RR}(X)$ is equipped with the cones
$C^k(X)$ (resp. $K^k(X)$) of strongly (resp. weakly) positive
elements as in Definition \ref{g1.5}.

Let now $X=\HH^n$ be the flat space. Consider the bundle
$\Ome^{k,k}(X):=X\times \Ome^{k,k}(\HH^n)$ over $X$, $0\leq k\leq
n$, where $\Ome^{k,k}(\HH^n)$ is as in Proposition
\ref{omega-def}. This bundle is isomorphic to the bundle
$\Lam^{2k,0}_{I,\RR}(X)$ via the isomorphism $\tau_k$ from Lemma
\ref{g5}. In \cite{alesker-valq} (see also the end of Section
\ref{linalg} of this article) we have described the cones of
strongly and weakly positive elements in the bundle
$\Ome^{k,k}(X)$. These cones correspond to the cones $C^k(X)$ and
$K^k(X)$ in $\lir(X)$ via the isomorphism $\tau_k$.

\section{Quaternionic plurisubharmonic functions \\ on
$\HH^n$.}\label{rempsh}

In this section we will remind the notion and basic facts on
quaternionic plurisubharmonic functions on $\HH^n$ following
\cite{alesker-bsm}, \cite{alesker-valq}.

Let $X$ be an open subset in $\HH ^n$.
\begin{definition}\label{def-psh}
A real valued function $u: X \to \RR$ is called quaternionic
plurisubharmonic if it is upper semi-continuous and its
restriction to any right {\itshape quaternionic} line is
subharmonic.
\end{definition}
 Recall that upper semi-continuity means that
 $u(x_0)\geq \underset{x\to x_0}{\limsup u(x)}$ for any $x_0\in X$.
 We will denote by $P(X)$ the class of plurisubharmonic
 functions in the open set $X$.

The class of all quaternionic plurisubharmonic functions in $X$
will be denoted by $P(X)$.

Sometimes we will abbreviate the term "quaternionic
plurisubharmonic" by just "plurisubharmonic".

\begin{proposition}[\cite{alesker-bsm}]\label{72}
Let $X\subset \HH^n$ be an open subset. Let $f\colon X\to \RR$ be
a $C^2$-smooth function. Then

(1) the matrix $\left(\frac{\6^2f}{\6\bar q_i\6  q_j}\right)$ is
hyperhermitian;

(2) the function $f$ is quaternionic plurisubharmonic if and only
if this matrix $\left(\frac{\6^2f}{\6\bar q_i\6  q_j}\right)$ is
non-negative definite.
\end{proposition}


Thus for any $C^2$-smooth function $f\colon X\to \RR$ the matrix
$\left(\frac{\6^2f}{\6\bar q_i\6  q_j}\right)$ is a continuous
section of the bundle $\Ome^{1,1}(X)$. Let us denote for brevity
$$D_2f:=\left(\frac{\6^2f}{\6\bar q_i\6 q_j}\right).$$
Thus by Proposition \ref{72} $f$ is quaternionic plurisubharmonic
function if and only if $D_2f$ take values in the cone of weakly
(= strongly) positive elements of $\Ome^{1,1}(X)$.

\begin{theorem}[\cite{alesker-valq}]\label{acln}
Let $X\subset \HH^n$ be an open subset. Let $0\leq k\leq n$. For
any functions $u_1,\dots,u_k\in C(X)\cap P(X)$ one can define a
$\Ome^{k,k}(X)\otimes |\ome_X|^*$-valued measure denoted by
$D_2u_1\cdot \dots D_2 u_k$ with values in the cone of strongly
positive elements which is uniquely characterized by the following
two properties:

1) if $u_1,\dots,u_k\in C^2(X)$ then the expression has the
obvious meaning;

2) if sequences $\{u_1^{(N)}\},\dots,\{u_k^{(N)}\}\subset C(X)\cap
P(X)$ are such that $u_i^{(N)}\to u_i\in C(X)\cap P(X)$ uniformly
on compact subsets of $X$ for any $i=1,\dots,k$ then
$$D_2u_1^{(N)}\cdot \dots D_2u_k^{(N)}
\overset{\mbox{weakly}}{\to}D_2u_1\cdot \dots D_2u_k$$ where the
convergence is understood is the sense of the weak convergence of
measures, i.e. in the space $(C_c(X,(\Ome^{k,k}(X))^*\otimes
|\ome_X|))^*$ equipped with the weak topology.
\end{theorem}

{\bf Remarks.} (a) It is easy to see that if $u_N\to u$ uniformly
on compact subsets, and $u_N\in C(X)\cap P(X)$ then $u\in C(X)\cap
P(X)$.

(b) Note that the real analogue of this result was proved by A.D.
Aleksandrov \cite{aleksandrov2}, and the complex analogue by
Chern, Levine, and Nirenberg \cite{chern-levine-nirenberg}.

(c) A special case of Theorem \ref{acln} with $k=n$ was proved
earlier by one of the authors \cite{alesker-bsm}.

\section{Plurisubharmonic functions on hypercomplex
manifolds.}\label{qpsh}
Let $(X^{4n},I,J,K)$ be a hypercomplex manifold.
\begin{definition}\label{2-1}
A continuous function $$h:X\to \RR$$ is called {\itshape
quaternionic plurisubharmonic} (or just psh) if $\6\6_J h$ is a
non-negative (generalized) section of $\Lambda^{2,0}_{I,\RR}(X)$
(where the non-negativity is understood in the sense of Section
\ref{vecbuncone}.
\end{definition}

Let us denote by $P'(X)$ the class of continuous psh functions on
$X$. Let us denote by $P''(X)$ the class of functions from $P'(X)$
with the following additional property: a function $h\in P'(X)$
belongs to $P''(X)$ if and only if any $x\in X$ has a neighborhood
$U\ni x$ and a sequence $\{h_N\}\subset P'(U)\cap C^2(U)$ such
that $h_N\overset{C^0(U)}{\to}h$. Thus $P''(X)\subset P'(X)$.

\begin{remark}\label{2-2}
If $X\subset \HH^n$ is an open subset then it is easy to see using
convolution with smooth non-negative functions that
$P''(X)=P'(X)=P(X)\cap C(X)$ where $P(X)$ is the class of
quaternionic plurisubharmonic functions as in Definition
\ref{def-psh}.
\end{remark}
{\bf Conjecture.} $P''(X)=P'(X)$ on any hypercomplex manifold $X$.

\begin{proposition}\label{2-3}
$P''(X)$ is closed under taking maximums.
\end{proposition}
To prove Proposition \ref{2-3} we will need a lemma.
\begin{lemma}\label{2-4}
Let $f,g \in P'(X)\cap C^2(X),\, f,g>0$. Then for any $1\leq
p<\infty$ $$(f^p+g^p)^\frac{1}{p}\in P'(X)\cap C^2(X).$$
\end{lemma}

Let us deduce Proposition \ref{2-3} from Lemma \ref{2-4}. Let $f,g
\in P''(X)$. Since the statement is local, adding a large constant
we may assume that $f,g>0$. Then there exist $\{f_N\},\,
\{g_N\}\subset P'(X)\cap C^2(X)$ such that
$f_N\overset{C^0(X)}{\to}f, \, g_N\overset{C^0(X)}{\to}g$.
Moreover we may assume that $f_N,g_N >0$ for all $N$. Then
$\max\{f_N,g_N\}\to \max\{f,g\}$. It remains to show that
$\max\{f_N,g_N\}\in P''(X)$. One has
$$\max\{f_N,g_N\}=\lim_{p\to \infty}(f^p+g^p)^\frac{1}{p}.$$
Hence Lemma \ref{2-4} implies Proposition \ref{2-3}.

{\bf Proof} of Lemma \ref{2-4}. Step 1. Let us reduce the
statement to the flat case $X=\HH^n$. The plurisubharmonicity of
$(f^p+g^p)^\frac{1}{p} $ for $f,g\in C^2(X)$ it is enough to check
pointwise. Let us fix $p\in X$.

The following lemma is well known (see e.g. \cite{salamon},
\cite{bryant}) though we will outline a proof for convenience of
the reader.
\begin{lemma}\label{2-5}
Let $x_0\in X$. There exists a neighborhood $U\ni x_0$ and a
diffeomorphism of $U$ onto an open subset $V\subset \HH^n$ such
that $x_0$ goes to 0 and
\begin{eqnarray*}
I(x)=I_0+O(|x|^2)\\
J(x)=J_0+O(|x|^2)\\
K(x)=K_0+O(|x|^2)
\end{eqnarray*}
where $I_0,J_0,K_0$ are the standard quaternionic structures on
$\HH^n$.
\end{lemma}
{\bf Proof.} The key point in the proof is the following result
due to Obata \cite{obata}: on the tangent bundle of hypercomplex
manifold $(X,I,J,K)$ there exists unique torsion free connection
$\nabla$ preserving the complex structure $I,J,K$:
$$\nabla I=\nabla J=\nabla K=0.$$ This connection is called the
Obata connection.

Let us consider the geodesic coordinates with respect to the Obata
connection $\nabla$ in a neighborhood of a point $x_0\in X$. Let
us denote by $I_0,J_0,K_0$ the flat complex structure in this
neighborhood. Thus
\begin{eqnarray}\label{2-5.1}
I(x_0)=I_0(x_0),\,J(x_0)=J_0(x_0), \,
K(x_0)=K_0(x_0).
\end{eqnarray}
Since the Obata connection is torsion free, the Cristoffel symbols
of $\nabla$ vanish at $x_0$. Hence
\begin{eqnarray}\label{2-5.2}
(\nabla I_0)(x_0)=(\nabla J_0)(x_0)=(\nabla K_0)(x_0)=0.
\end{eqnarray}
Then (\ref{2-5.1}), (\ref{2-5.2}) imply the proposition. \qed

It is enough to know $f,g$ in the 2-jet neighborhood of $x_0$, and
the complex structures in 1-jet neighborhood of $x_0$. Hence by
Lemma \ref{2-5} we reduce Lemma \ref{2-4} to the flat case.

Step 2. Assume that $X=\HH^n$. Let us fix an arbitrary right
$\HH$-line $L$. We have to check that $(f^p+g^p)^\frac{1}{p}|_L$
is subharmonic. Hence it remains to prove the following lemma.
\begin{lemma}\label{2-6}
Let $f,g>0$ be continuous subharmonic functions on a Euclidean
space $L$. Then for any $1\leq p<\infty$ the function
$(f^p+g^p)^\frac{1}{p}$ is subharmonic.
\end{lemma}

To prove Lemma \ref{2-6} (which is well known) let us fix a sphere
$S$ with a center $x_0\in L$. Let $F(x):=(f(x),g(x))\in \RR^2$. We
have
\begin{eqnarray*}
(f^p(x_0)+g^p(x_0))^\frac{1}{p}\leq \left((\int_S f(x))^p+(\int _S
g(x))^p\right)^\frac{1}{p}=\\
||\int_S F(x)||_{l^p}\leq
\int_S||F(x)||_{l^p}=\int_S(f(x)^p+g(x)^p)^\frac{1}{p}.
\end{eqnarray*}
Hence $(f^p+g^p)^\frac{1}{p}$ is subharmonic. Lemma \ref{2-6} is
proved. Hence Lemma \ref{2-4} is proved too. \qed

\begin{lemma}\label{2-7}
Let $X$ be a hypercomplex manifold, $\dim _\RR X=4n$. For any
$\xi\in C^\infty_0(X,\Omega^k),\, \eta\in
C^\infty_0(X,\Omega^{4n-k-1})$  one has
$$\int_X\xi\wedge \6_J\eta =-\int_X \6_J\xi \wedge \eta.$$
\end{lemma}
{\bf Proof.} First observe that for any top degree form $\rho\in
\Omega^{4n}(X)$ one has
$$J\circ \rho=\rho.$$
Next for any forms $\alpha,\beta$ one has
$$J\circ(\alpha \wedge\beta)=(J\circ \alpha)\wedge(J\circ
\beta).$$ Using this and the Stokes formula we get
\begin{eqnarray*}
\int \xi\wedge \6_J\eta=\int\xi \wedge (J^{-1}\circ \bar\6\circ
J)\eta=\\
\int J\xi \wedge \bar\6\circ J\circ \eta=-\int \bar\6\circ
J\circ\xi \wedge J\eta=\\
-\int J^{-1}\circ \bar\6\circ J\circ \xi \wedge \eta=-\int \6_J\xi
\wedge \eta.
\end{eqnarray*}

\qed

\begin{proposition}\label{2-8}
Let $\{h_N\}\subset C(X)$. Let $h_N\overset{C^0}{\to}h$. Then

(1) if for any $N$ $h_N\in P'(X)$ then $h\in P'(X)$;

(2) if for any $N$ $h_N\in P''(X)$ then $h\in P''(X)$.
\end{proposition}
{\bf Proof.} Part (2) easily follows from part (1). Thus let us
prove part (1). We have to check that $h$ is psh. Let $\phi\in
C^\infty_0(X,S_\HH ^*\otimes \ox) \simeq
C^\infty(X,\Lambda^{2n-2,2n}_{I,\RR})$ be a non-negative section.
Then by Lemma \ref{2-7} one gets
\begin{eqnarray*}
\int \6\6_J h\wedge \phi =\int h\6\6_J \phi=\\
\lim_{N\to \infty}\int h_N\6\6_J \phi=\lim_{N\to \infty} \int
\6\6_J h_N \wedge \phi \geq 0.
\end{eqnarray*}
The following result is an analogue of  the theorems of
Aleksandrov \cite{aleksandrov2} and Chern-Levine-Nirenberg
\cite{chern-levine-nirenberg}.
\begin{theorem}\label{2-9}
Let $0<k\leq n$. For any $h_1,\dots,h_k\in P''(X)$ one can define
a non-negative generalized section of $\Lambda^{2k}_{I,\RR}(X)$
denoted by $\6\6_J h_1\wedge \dots\wedge\6\6_J h_k$ which is
uniquely characterized by the following two properties:

(1) if $h_1,\dots,h_k\in C^2(X)$ then the definition is clear;

(2) if $\{h_i^N\}\subset C^2(X)$, $h_i^N\overset{C^0}{\to}h_i$ as
$N\to \infty$, $i=1,\dots, k$ then $h_i\in P''(X)$ and
$$\6\6_J h_1^N\wedge \dots\wedge \6\6_J h_k^N\to \6\6_J h_1\wedge \dots\wedge \6\6_J h_k$$
in the weak topology on $\Lambda^{2k}_{I,\RR}(X)\otimes
|\ome_X|^*)$-valued measures.
\end{theorem}
\begin{remark}
Theorem \ref{acln} is a special case of Theorem \ref{2-9}.
\end{remark}
To prove Theorem \ref{2-9} we will need a lemma.
\begin{lemma}\label{2-10}
Let $0\leq k\leq n$. Fix a compact subset $K\subset X$ and its
compact neighborhood $\tilde K\supset K$. Then there exists a
constant $C$ such that for any $h_1,\dots,h_k\in P'(X)\cap C^2(X)$
one has
$$||\6\6_J h_1\wedge \dots \wedge \6\6_J h_k||_{L^1(K)}\leq C
\prod_{i=1}^k||h_i||_{C(\tilde K)}.$$
\end{lemma}
 {\bf Proof} of Lemma \ref{2-10}. Let us prove the lemma by
induction in $k$. For $k=0$ the lemma is trivial. Let us assume
that the lemma is true for $k-1$ and let us prove it for $k$. Let
us fix a compact neighborhood $\hat K$ of $K$ such that $\tilde K$
is a neighborhood of $\hat K$:
$$K\subset \hat K\subset \tilde K.$$
 Let us fix $\gamma \in
C^\infty_0(X,\Lambda^{2n-2k,2n}_I(X))$ which takes values in the
cone of weakly positive elements, and moreover the restriction of
$\gamma$ to $K$ takes values in the interior of this cone, and
$\gamma|_{X\backslash \hat K}\equiv 0$. Then
\begin{eqnarray*}
||\6\6_J h_1\wedge \dots \wedge \6\6_J h_k||_{L^1(K)}\leq
C_1\int_X \6\6_J h_1\wedge \dots \wedge \6\6_J h_k\wedge
\gamma=\\
C_1 \int_X h_k\6\6_J (\6\6_J h_1\wedge \dots \wedge \6\6_J h_{k-1}\wedge \gamma)=\\
C_1\int_{\hat K}h_k\6\6_J h_1\wedge \dots \wedge \6\6_J h_{k-1}\wedge \6\6_J \gamma \leq\\
C_2||h_k||_{C(\hat K)} \int _{\hat K}|\6\6_J h_1\wedge \dots
\wedge \6\6_J h_{k-1}\wedge \6\6_J \gamma|\leq\\
C_3||h_k||_{C(\hat K)}||\6\6_J h_1\wedge \dots \wedge \6\6_J
h_{k-1}||_{L^1(\hat K)}\leq C\prod_{i=1}^k||h_i||_{C(\tilde K)}
\end{eqnarray*}
where the last inequality follows from the assumption of
induction. Lemma \ref{2-10} is proved. \qed

{\bf Proof} of Theorem \ref{2-9}. Let $h_1,\dots,h_k\in P''(X)$.
Let us choose \[
  \{h_i^N\}_{N=1}^\infty \subset P'(X)\cap C^2(X),\, i=1,\dots, k
\] such that $h_i^N \overset{C^0}{\to}h_i$ for each
$i=1,\dots, k$. Let us show that $\prod_{i=1}^k\6\6_J h_i^N$
converges weakly to a strongly non-negative measure with values in
$\Lambda^{2k,0}_I(X)$. Since the $\Lambda^{2k,0}_I(X)\otimes
|\ome_X|^* $-valued measures $\prod_{i=1}^k \6\6_J h^N_i$ are
non-negative and locally bounded (by Lemma \ref{2-10}) there
exists a subsequence $\{N_l\}$ such that $\prod_{i=1}^k \6\6_J
h^{N_l}_i$ converges weakly to a non-negative
$\Lambda^{2k,0}_I(X)\otimes |\ome_X|^*$-valued measure $\mu$. Let
us show that $\mu$ does not depend on a choice of a convergent
subsequence. Let us show it by induction in $k$. For $k=0$ this is
trivial. Let us assume that the statement is true for $k-1$ and
let us prove it for $k$. Let $\nu$ be another weak limit of some
subsequence of the sequence $\{\prod_{i=1}^k \6\6_J h^N_i\}$. It
is enough to check that for any $\phi\in
C^\infty_0(X,\Lambda^{2n-2k,2n}(X))$
\begin{eqnarray}\label{2-11}
\int_X\mu\wedge \phi= \int_X\nu\wedge \phi.
\end{eqnarray}
We have
\begin{eqnarray}\label{2-11.5}
\int\mu\wedge\phi=\lim_{l\to \infty} \int\prod_{i=1}^k \6\6_J
h^{N_l}_i\wedge \phi=\\\label{2-11.55} \lim_{l\to \infty} \int
h_k^{N_l} \cdot (\prod_{i=1}^{k-1} \6\6_J h^{N_l}_i)\wedge \6\6_J
\phi.
\end{eqnarray}
By the assumption of induction the sequence
$\{g_N:=\prod_{i=1}^{k-1} \6\6_J h^{N_l}_i\}$ is weakly
convergent. Let us denote its weak limit by $g$. The sequence
$\{f_N:=h_k^N\cdot \6\6_J \phi\}$ has uniformly bounded support
and converges uniformly (i.e. in $C^0$-topology) to $f:=h_k\cdot
\6\6_J \phi$. Thus the existence of the limit in (\ref{2-11.55})
follows from the following known lemma.
\begin{lemma}\label{2-11.6}
Let $M$ be a compact topological space. Let $E\to M$ be a vector
bundle. Let $\{f_N\}\subset C(M,E)$ be a sequence of sections
which converges to $f\in C(M,E)$ uniformly on $M$. Let
$\{g_N\}\subset C(M,E)^*$ be a sequence in the dual space which is
weakly convergent to $g\in C(M,E)^*$. Then
$$g_N(f_N)\to g(f).$$
\end{lemma}
Let us postpone the proof of this lemma which is well known. This
lemma implies that $\lim_{N\to \infty}\int h_k^{N} \cdot
(\prod_{i=1}^{k-1} \6\6_J h^{N}_i)\wedge \6\6_J \phi$ does exist,
and by the same argument it should be equal to $\int\nu\wedge
\phi$. Hence the equality (\ref{2-11}) follows. Hence there exists
a weak limit of the sequence $\prod_{i=1}^k\6\6_J h_i^N$.

It remains to show that if $h_1,\dots,h_k\in P'(X)\cap C^2(X)$
then the limit is equal to $\prod_{i=1}^k\6\6_J h_i$. Let us show
this by induction in $k$. For $k=0$ this is trivial. Let us make
the induction step. Let us denote by $\mu$ the weak limit of
$\prod_{i=1}^k\6\6_J h_i^N$. Let us fix $\phi\in
C^\infty_0(X,\Lambda^{2n-2k,2n}(X))$. We have
\begin{eqnarray*}
\int\mu\wedge\phi=\lim_{N\to \infty}
\int\prod_{i=1}^k\6\6_J h_i^N\wedge \phi=\\
\lim_{N\to \infty} \int h_k^N \cdot\prod_{i=1}^{k-1}\6\6_J
h_i^N\wedge\6\6_J \phi.
\end{eqnarray*}
By Lemma \ref{2-11.6} the last limit is equal to
\begin{eqnarray}\label{zzz}
\int h_k\cdot\prod_{i=1}^{k-1}\6\6_J h_i\wedge\6\6_J \phi.
\end{eqnarray}
Let us show that (\ref{zzz}) is equal to $\prod_{i=1}^k\6\6_Jh_i
\wedge \phi$. Using approximation and Lemma \ref{2-11.6}, we may
assume that $h_1,\dots,h_k\in P'(X)\cap C^2(X)$. This case follows
from Lemma \ref{2-7}. Theorem \ref{2-9} is proved. \qed

{\bf Proof} of Lemma \ref{2-11.6}. By the Banach-Steinhauss
theorem the sequence $\{g_N\}$ in bounded in the $C(M,E)^*$-norm.
Then we have
\begin{eqnarray*}
|g_N(f_N)-g(f)|\leq |g_N(f_N-f)|+|(g_N-g)(f)|\leq\\
C||f_N-f||_{C(M,E)}+|(g_N-g)(f)|\to 0.
\end{eqnarray*}
\qed

\section{Quaternionic plurisubharmonic functions \\ and
HKT-geometry.}\label{dhkt}

In this section we present a geometric interpretation of
quaternionic strictly plurisubharmonic functions on hypercomplex
manifolds as local potentials of HKT-metrics (see Definition
\ref{def-hkt-metr}). This result is analogous to the classical
well known fact that the complex strictly plurisubharmonic
functions on complex manifolds are precisely local potentials of
K\"ahler metrics.

Let $(X^{4n},I,J,K)$ be a hypercomplex manifold. Remind that we
have the natural identification
$$t\colon \Lam^{2,0}_{I,\RR}(X)\tilde \to S_\HH(X)$$
defined by (\ref{n7}) in Section \ref{operators}. The following
proposition which we call a local $\6\6_J$-lemma, is a rather
straightforward application of the main result of
\cite{banos-swann} (see also \cite{baston}).


\begin{proposition}\label{ddbar}
Let $\Ome\in C^\infty(X,\Lam^{2,0}_{I,\RR})$. Then locally on $X$
the form $\Ome$ can be presented in a form
$$\Ome=\6\6_Jf$$
with $f$ being a $C^\infty$-smooth real valued function if and
only if $\6\Ome=0$.
\end{proposition}
{\bf Proof.} The only if part is clear:
$$\6(\6\6_J)\Ome)=\6^2(\6_J\Ome)=0.$$
Let us prove the converse statement. We will use heavily the
result of \cite{banos-swann} (which is based in turn on
\cite{mamone-salamon}, see also \cite{baston}). The article
\cite{banos-swann} uses a convention that the complex structures
$I,J,K$ act on the {\itshape left} on the tangent vectors to $X$.
Though in this article we use the opposite convention, in the
proof of this proposition we will change our convention in order
to cite the results of \cite{banos-swann} without change of signs
and normalizations.

For a complex structure $\ci$ compatible with the quaternionic
structure (namely $\ci$ has a form $\ci=aI+bJ+cK,\, a,b,c\in \RR,
a^2+b^2+c^2=1$) let us define the action of $\ci$ on $k$ forms
$$(\ci\circ\ome)(X_1,\dots,X_k)=\ome(\ci^{-1} X_1,\dots,\ci^{-1} X_k).$$
Define differentials on $k$-forms
\begin{eqnarray*}
d_\ci\ome:=-\ci^{-1}\circ d\circ \ci\circ\ome.
\end{eqnarray*}
Note that the operators $d_\ci$ coincide with those considered in
Section 2 of \cite{banos-swann}. It is known (and one can be
easily check it by a straightforward computation) that
$d,d_I,d_J,d_K$ anticommute. Also one has
$$\6=\frac{1}{2}(d+\sqrt{-1}d_I),\, \bar \6=\frac{1}{2}(d-\sqrt{-1}d_I).$$
Next one obtains
\begin{eqnarray*}
\6\6_J=\frac{1}{4}(d+\sqrt{-1}d_I)\circ
J^{-1}\circ(d-\sqrt{-1}d_I)\circ
J=\\\frac{1}{4}(d+\sqrt{-1}d_I)(J^{-1}\circ d\circ
J+\sqrt{-1}J^{-1}I^{-1}\circ d\circ IJ)=\\
\frac{1}{4}(d+\sqrt{-1}d_I)(-d_J-\sqrt{-1}d_K)=\\
-\frac{1}{4}((dd_J-d_Id_K)+\sqrt{-1}(dd_K+d_Id_J))=\\
-\frac{1}{4}\left((dd_J+d_Kd_I)+\sqrt{-1}(dd_K+d_Id_J)\right).
\end{eqnarray*}
Thus we get
\begin{eqnarray}\label{ot}
\6\6_J=-\frac{1}{4}\left((dd_J+d_Kd_I)+\sqrt{-1}(dd_K+d_Id_J)\right).
\end{eqnarray}

Let us denote by $g:=t(\Ome)$ the smooth section of
$\Lam^{2,0}_{I,\RR}(X)$ corresponding to $\Ome$. For a complex
structure $\ci$ let us define a 2-form
$$F_\ci(X,Y)=g(\ci X,Y).$$ Then
\begin{eqnarray}\label{ow}
\Ome=c\cdot (F_J+\sqrt{-1}F_K)
\end{eqnarray} where $c$ is a
non-zero normalizing constant which we will not write down
explicitly. (Note also that in the right hand side of (\ref{ow})
we have the sign "plus" instead of "minus" as previously since now
we work with left hypercomplex structures in the opposite to our
previous conventions).

By the main theorem of \cite{banos-swann} locally on $X$ there
exists an infinitely smooth function $\mu$ such that
$$F_I=\frac{1}{2}(dd_I+d_Jd_K)\mu.$$
By Remark at the beginning part of Section 2 of \cite{banos-swann}
the last identity is equivalent to each of the following two
identities:
\begin{eqnarray*}
F_J=\frac{1}{2}(dd_J+d_Kd_I)\mu;\\
F_K=\frac{1}{2}(dd_K+d_Id_J)\mu.
\end{eqnarray*}
Substituting the last two identities into (\ref{ow}) we obtain
$$\Ome=\frac{c}{2}((dd_J+d_Kd_I)+\sqrt{-1}(dd_K+d_Id_J))\mu.$$
Using (\ref{ot}) we deduce
$$\Ome=c'\cdot \6\6_J\mu.$$
Thus Proposition \ref{ddbar} is proved. \qed


\begin{proposition}
(1) Let $f$ be an infinitely smooth strictly plurisubharmonic
function on a hypercomplex manifold $(X,I,J,K)$. Then $t(\6\6_J
f)$ is an HKT-metric.

(2) Conversely assume that $g$ is an HKT-metric. Then any point
$x\in X$ has a neighborhood $U$ and an infinitely smooth strictly
plurisubharmonic function $f$ on $U$ such that $g=t(\6\6_J f)$ in
$U$.
\end{proposition}
{\bf Proof.} This is an immediate consequence of Proposition
\ref{ddbar} and the definition of plurisubharmonic function. \qed

\hfill

Non-negativity of a (2,0)-form $\rho \in \Lambda^{2,0}_{I,
\RR}(X)$ was defined and explored at great length in
\cite{_V:reflexive_}, under the name ``$K$-positivity''. There it
was related to a notion of stability for coherent sheaves on
hyperk\"ahler manifolds, and used to establish stability in some
important cases. Quaternionic analogues of some fundamental
results on positive currents were proven; in particular, a
quaternionic version of Sibony's lemma on extensions
(\cite{_Sibony_}) and a version of Skoda-El Mir theorem.

\hfill

\begin{theorem}\label{_sibony_el_mir_Theorem_}
Let $(X,I,J,K)$ be a hyperk\"ahler manifold, $Z\subset (X,K)$ a
complex subvariety, $codim_{\Bbb C}(Z, X)\leq 3$, and $\rho \in
\Lambda^{2,0}_{I, \RR}(X)$ a $\6$-closed non-negative form on
$X\backslash Z$. Then $\rho$ is locally $L^1$-integrable on $X$,
and the corresponding current on $X$ is $\6$-closed.
\end{theorem}

{\bf Proof:} This is Proposition 7.5 from \cite{_V:reflexive_}.
\qed

\hfill

\vskip 0.7cm

{\small

\hfill

\noindent {\sc Semyon Alesker\\
{ \normalsize Department of Mathematics, Tel Aviv University,
Ramat Aviv}
 \\  { \normalsize 69978 Tel Aviv,
Israel }
\\ \tt semyon@post.tau.ac.il\\

\noindent \sc Misha Verbitsky\\
University of Glasgow, Department of Mathematics, \\
15 University Gardens, Glasgow G12 8QW, Scotland.}\\
\   \\
{\sc  Institute of Theoretical and
Experimental Physics \\
B. Cheremushkinskaya, 25, Moscow, 117259, Russia }\\
\  \\
\tt verbit@maths.gla.ac.uk, \ \  verbit@mccme.ru
}


\begin{thebibliography}{99}

\bibitem[Al]{aleksandrov2}
 Aleksandrov, A. D.; Dirichlet's problem for the equation
 $ Det\,||z_{ij}|| =\varphi (z_{1},\cdots,z_{n},z, x_{1},\cdots, x_{n})$. I. (Russian) Vestnik
Leningrad. Univ. Ser. Mat. Meh. Astr. 13 1958 no. 1, 5--24.

\bibitem[A1]{alesker-bsm}
Alesker, Semyon; Non-commutative linear algebra and
plurisubharmonic functions of quaternionic variables. Bull. Sci.
Math., 127 (2003), no. 1, 1--35. also: math.CV/0104209

\bibitem[A2]{alesker-ma}
Alesker, Semyon; Quaternionic Monge-Ampère equations. J. Geom.
Anal. 13 (2003), no. 2, 205--238. also: math.CV/0208005.

\bibitem[A3]{alesker-valq}
Alesker, Semyon; Valuations on convex sets, non-commutative
determinants, and pluripotential theory. Adv. Math. 195 (2005),
561-595. also: math.MG/0401219.

\bibitem[As]{aslaksen}
Aslaksen, Helmer;
 Quaternionic determinants. Math. Intelligencer 18 (1996), no. 3, 57--65.

\bibitem[B]{baston}
Baston, R. J.; Quaternionic complexes. J. Geom. Phys. 8 (1992),
no. 1-4, 29--52.

\bibitem[BS]{banos-swann}  Banos, Bertrand; Swann, Andrew;
 Potentials for hyper-K\"ahler metrics with torsion.
 Classical Quantum Gravity 21 (2004), no. 13, 3127--3135.
\bibitem[Bo]{boyer}
Boyer, Charles P.; A note on hyper-Hermitian four-manifolds. Proc.
Amer. Math. Soc. 102 (1988), no. 1, 157--164.

\bibitem[Br]{bryant} Bryant, Robert;
Classical, exceptional, and exotic holonomies: a status report, in
Actes de la Table Ronde de G\'eom\'etrie Diff\'erentielle (Luminy,
1992), S\'emin. Congr., vol. 1 (1996), pp. 93-165, Soc. Math.
France, Paris.

\bibitem[CLL]{chern-levine-nirenberg}
 Chern, Shiing-shen; Levine, Harold I.; Nirenberg, Louis;
 Intrinsic norms on a complex manifold.
 1969 Global Analysis (Papers in Honor of K. Kodaira) pp. 119--139
Univ. Tokyo Press, Tokyo.

\bibitem[GP]{_Gra_Poon_}
Grantcharov, Gueo; Poon, Yat Sun; Geometry of hyper-K\"ahler
connections with torsion, Comm. Math. Phys. 213 (2000), no. 1,
19--37.


\bibitem[GV]{generalized-functions}
Gelfand, I. M.; Vilenkin, N. Ya.; Generalized functions. Vol. 4.
Applications of harmonic analysis. Translated from the Russian by
Amiel Feinstein. Academic Press [Harcourt Brace Jovanovich,
Publishers], New York-London, 1964 [1977].

\bibitem[GRW]{gelfand-retakh-wilson}
Gelfand, Israel; Retakh, Vladimir; Wilson, Robert Lee;
Quaternionic quasideterminants and determinants. Lie groups and
symmetric spaces, 111--123, Amer. Math. Soc. Transl. Ser. 2, 210,
Amer. Math. Soc., Providence, RI, 2003.

\bibitem[Ha]{harvey}
Harvey, Reese; Holomorphic chains and their boundaries. Several
complex variables (Proc. Sympos. Pure Math., Vol. XXX, Part 1,
Williams Coll., Williamstown, Mass., 1975), pp. 309--382. Amer.
Math. Soc., Providence, R. I., 1977.

\bibitem[He]{henkin}
Henkin, Gennadii; Private communication.

\bibitem[Ho]{howe}
Howe, Roger; Remarks on classical invariant theory. Trans. Amer.
Math. Soc. 313 (1989), no. 2, 539--570.

\bibitem[HP]{howe-papa}
Howe, P. S.; Papadopoulos, G.; Twistor spaces for hyper-K\"ahler
manifolds with torsion. Phys. Lett. B 379 (1996), no. 1-4, 80--86.

\bibitem[L]{lelong}
Lelong, Pierre; Fonctions plurisousharmoniques et formes
diff\'erentielles positives. (French) Gordon \& Breach,
Paris-London-New York (Distributed by Dunod \'editeur, Paris)
1968.
\bibitem[MCSa]{mamone-salamon}
Mamone Capria, M.; Salamon, S. M.; Yang-Mills fields on
quaternionic spaces. Nonlinearity 1 (1988), no. 4, 517--530.

\bibitem[M]{moore}
Moore,E.H.; On the determinant of an hermitian matrix of
quaternionic elements. Bull. Amer. Math. Soc. 28 (1922), 161-162
\bibitem[Ob]{obata}
Obata, Morio; Affine connections on manifolds with almost complex,
quaternion or Hermitian structure.
 Jap. J. Math. 26, 1956 43--77.
\bibitem[Sa]{salamon}
Salamon, Simon M.; Differential geometry of quaternionic
manifolds. Annales Scientifiques de l'École Normale Supérieure
Sér. 4, 19 no. 1 (1986), p. 31-55.

\bibitem[Sib]{_Sibony_}
Sibony, Nessim; Quelques problemes de prolongement de courants en
analyse complexe, Duke Math. J. 52, 157-197 (1985).



\bibitem[V1]{_V:reflexive_}
Verbitsky, Misha; Hyperholomorpic connections   on coherent
sheaves and stability, 40 pages, math.AG/0107182


\bibitem[V2]{verbitsky-hkt}
Verbitsky, Misha; HyperK\"ahler manifolds with torsion,
supersymmetry and Hodge theory. Asian J. Math. 6 (2002), no. 4,
679--712.



\end{thebibliography}
\end{document}